\documentclass[11pt]{article}
\usepackage{amsfonts, amsthm, amssymb, eucal, amsbsy}
\usepackage{amsmath, chatterjee}

\begin{document}
\newcommand{\tx}{\tilde{X}}
\newcommand{\ty}{\tilde{Y}}

\title{An error bound in the Sudakov-Fernique inequality}
\author{Sourav Chatterjee}
\maketitle

\begin{abstract}
We obtain an asymptotically sharp error bound in the classical Sudakov-Fernique comparison inequality for finite collections of gaussian random variables. Our proof is short and self-contained, and gives an easy alternative argument for the classical inequality, extended to the case of non-centered processes. 
\end{abstract}

\section{Statement of the result}
Gaussian comparison inequalities are among the most important tools in the theory of gaussian processes, and the Sudakov-Fernique inequality (named after Sudakov \cite{sudakov71, sudakov76} and Fernique \cite{fernique75}) is perhaps the most widely used member of that class.

We will concentrate on the Sudakov-Fernique inequality in this article; general discussions about comparison inequalities can be found in Adler \cite{adler90}, Fernique \cite{fernique97}, Ledoux \& Talagrand \cite{ledoux85}, and Lifshits~\cite{lifshits95}.

The classical Sudakov-Fernique inequality goes as follows:
\begin{thm}\label{sf}
\textup{[Sudakov-Fernique inequality]} Let $\{X_i, i\in I\}$ and $\{Y_i, i\in I\}$ be
two centered gaussian processes indexed by the same indexing set $I$. Suppose that both the processes are almost surely bounded. For
each $i,j \in I$, let $\gamma^X_{ij}= \ee(X_i - X_j)^2$ and
$\gamma^Y_{ij} = \ee(Y_i - Y_j)^2$. If $\gamma^X_{ij} \le \gamma^Y_{ij}$ for all $i,j$, then $\ee(\sup_{i\in I} X_i) \le \ee(\sup_{i\in I} Y_i)$. 
\end{thm}
As mentioned before, this inequality is attributed to Sudakov \cite{sudakov71, sudakov76} and Fernique \cite{fernique75}. Later proofs were given in Alexander \cite{alexander85} and an unpublished work of S.\ Chevet. Important variants were proved by Gordon \cite{gordon85, gordon87, gordon92} and Kahane \cite{kahane86}.
More recently, Vitale \cite{vitale00} has shown, through a clever argument, that we only need $\ee(X_i)=\ee(Y_i)$ instead of $\ee(X_i)=\ee(Y_i)=0$ in the hypothesis of Theorem \ref{sf}.
We will prove the following result, which gives an sharp error bound when the indexing set is finite, and also contains Vitale's extension of the Sudakov-Fernique inequality.
\begin{thm}\label{max}
Let $(X_1,\ldots,X_n)$ and $(Y_1,\ldots, Y_n)$ be
gaussian random vectors with $\ee(X_i) =\ee(Y_i)$ for each $i$. For
$1\le i, j\le n$, let $\gamma^X_{ij}= \ee(X_i - X_j)^2$ and
$\gamma^Y_{ij} = \ee(Y_i - Y_j)^2$, and let $\gamma = \max_{1\le
  i,j\le n} |\gamma^X_{ij} - \gamma^Y_{ij}|$. Then
\[
|\ee(\max_{1\le i\le n} X_i) - \ee(\max_{1\le i\le n} Y_i)| \le
\sqrt{\gamma \log n}.
\]
Moreover, if $\gamma^X_{ij} \le \gamma^Y_{ij}$ for all $i,j$, then
$\ee(\max_i X_i) \le \ee(\max_i Y_i)$. 
\end{thm}
The asymptotic sharpness of the error bound is easy to see from the case where all the $X_i$'s are independent standard normals and all the $Y_i$'s are zero.

\section{Proof}
We first need to state the following well-known 
``integration by parts'' lemma:
\begin{lmm}\label{stein}
If $F:\rr^n \ra \rr$ is a $C^1$ function of moderate growth at
infinity, and $\bbx=(X_1,\ldots,X_n)$ is a centered Gaussian random
vector, then for any $1\le i \le n$,
\[
\ee(X_i F(\bbx)) = \sum_{j=1}^n \ee(X_iX_j) \ee\biggl(\frac{\partial
  F}{\partial x_i} (\bbx)\biggr).
\]
\end{lmm}
A proof of this lemma can be found in the appendix of
\cite{talagrand03}, for example. \\
\\ 
\indent {\sc Proof of Theorem \ref{max}.} Let $\bbx = (X_1,\ldots,X_n)$ and $\bby = (Y_1,\ldots,Y_n)$. Without loss of generality, we may
assume that $\bbx$ and  
$\bby$ are defined on the same probability space and are
independent. Fix $\beta >0$, and define $F_\beta : \rr^n \ra \rr $ as: 
\[
F_\beta(\bx) := \beta^{-1}\log\biggl(\sum_{i=1}^n e^{\beta
  x_i}\biggr).
\]
(Note that $\bx$ denotes the vector $(x_1,\ldots,x_n)$, a convention
that we shall follow throughout.) 
Now, for each $i$, let $\mu_i = \ee(X_i) = \ee(Y_i)$, $\tx_i = X_i - \mu_i$, and
$\ty_i = Y_i - \mu_i$. For $1\le i,j\le n$, let $\sigma^X_{ij} =
\ee(\tx_i\tx_j)$ and $\sigma^Y_{ij} = \ee(\ty_i\ty_j)$. For $0\le 
t\le 1$ define the random vector $\bbz_t = (Z_{t,1},\ldots, Z_{t,n})$
as
\[
Z_{t,i} = \sst \tx_i + \st \ty_i + \mu_i.
\]
For all $t\in [0,1]$, let $\varphi(t) =
\ee(F_\beta(\bbz_t))$. Then $\varphi$ is differentiable, and
\begin{align*}
\varphi^\prime(t) &= \ee \biggl[\sum_{i=1}^n \frac{\partial
  F_\beta}{\partial x_i} (\bbz_t) \biggl(\frac{\ty_i}{2\st} -
\frac{\tx_i}{2\sst}\biggr) \biggr].
\end{align*}
Again, for any $i$, Lemma \ref{stein} gives us
\[
\ee\biggl(\frac{\partial  F_\beta}{\partial x_i} (\bbz_t) \tx_i\biggr) 
= \sst \sum_{j=1}^n \sigma^X_{ij} \ee \biggl(\frac{\partial^2
  F_\beta}{\partial x_j \partial x_i} (\bbz_t) \biggr)
\]
and
\[
\ee\biggl(\frac{\partial  F_\beta}{\partial x_i} (\bbz_t) \ty_i\biggr) 
= \st \sum_{j=1}^n \sigma^Y_{ij} \ee \biggl(\frac{\partial^2
  F_\beta}{\partial x_j \partial x_i} (\bbz_t) \biggr).
\]
Combining, we have
\[
\varphi^\prime(t) = \frac{1}{2}\sum_{1\le i, j\le n} \ee
\biggl(\frac{\partial^2   F_\beta}{\partial x_j \partial x_i} (\bbz_t)
\biggr) (\sigma^Y_{ij} - \sigma^X_{ij}).
\]
Now 
\[
\frac{\partial F_\beta}{\partial x_i}(\bx) = p_i(\bx) :=
\frac{e^{\beta x_i}}{\sum_{j=1}^n e^{\beta x_j}}.
\]
Note that for each $\bx \in \rr^n$, the numbers $p_1(\bx),\ldots
p_n(\bx)$ as defined above are nonnegative and sum to $1$. In other
words, they induce a probability measure on $\{1,2,\ldots,n\}$. It is straightforward to verify that 
\[
\frac{\partial^2 F_\beta}{\partial x_j \partial x_i}(\bx) = 
\begin{cases}
\beta(p_i(\bx) - p_i(\bx)^2) & \text{ if } i=j,\\
-\beta p_i(\bx)p_j(\bx) & \text{ if } i\ne j.
\end{cases}
\]
Thus,
\begin{align*}
& \sum_{1\le i, j\le n} \frac{\partial^2   F_\beta}{\partial
  x_j \partial x_i} (\bx) (\sigma^Y_{ij} -
\sigma^X_{ij})\\
&= \beta \sum_{i=1}^n p_i(\bx)(\sigma^Y_{ii} - \sigma^X_{ii}) 
- \beta\sum_{1\le i,j\le n} p_i(\bx)p_j(\bx)(\sigma^Y_{ij} - \sigma^X_{ij}).
\end{align*}
Now observe that since $\sum_{i=1}^n p_i(\bx) = 1$, therefore
\[
\sum_{i=1}^n p_i(\bx)(\sigma^Y_{ii} - \sigma^X_{ii}) =
\frac{1}{2}\sum_{1\le   i,j\le n} p_i(\bx)p_j(\bx)(\sigma^Y_{ii}
- \sigma^X_{ii} + \sigma^Y_{jj} - \sigma^X_{jj}).
\]
Combining, we have
\begin{align*}
& \sum_{1\le i, j\le n} \frac{\partial^2   F_\beta}{\partial
  x_j \partial x_i} (\bx)  (\sigma^Y_{ij} -
\sigma^X_{ij})\\
&= \frac{\beta}{2} \sum_{1\le i, j \le n} p_i(\bx)
p_j(\bx)\bigl[(\sigma^Y_{ii} + \sigma^Y_{jj} -
2\sigma^Y_{ij}) - (\sigma^X_{ii} + \sigma^X_{jj} -
2\sigma^X_{ij})\bigr].
\end{align*}
Now note that
\[
\sigma^X_{ii} + \sigma^X_{jj} - 2\sigma^X_{ij} =
\ee(\tx_i - \tx_j)^2 = \ee(X_i - X_j)^2 - (\mu_i - \mu_j)^2
\]
and similarly
\[
\sigma^Y_{ii} + \sigma^Y_{jj} - 2\sigma^Y_{ij} =
\ee(\ty_i - \ty_j)^2 = \ee(Y_i - Y_j)^2 - (\mu_i - \mu_j)^2.
\]
Therefore,
\begin{align*}
\sum_{1\le i, j\le n} \frac{\partial^2   F_\beta}{\partial
  x_j \partial x_i} (\bx) (\sigma^Y_{ij} -
\sigma^X_{ij})&= \frac{\beta}{2} \sum_{1\le i, j \le n} p_i(\bx)
p_j(\bx)(\gamma^Y_{ij} - \gamma^X_{ij}).
\end{align*}
Thus, if $\gamma^X_{ij} \le\gamma^Y_{ij}$ for all $i,j$, then
$\varphi^\prime(t) \ge 0$ for each $t$, which implies
\begin{equation}\label{ineq1}
\ee(F_\beta(\bby)) = \varphi(1) \ge \varphi(0) = \ee(F_\beta(\bbx)).
\end{equation}
Now observe that
\begin{align}
\max_i x_i &= \beta^{-1} \log e^{\beta \max_i x_i} \nonumber\\
&\le \beta^{-1} \log \biggl(\sum_i e^{\beta x_i}\biggr) \nonumber\\
&\le \beta^{-1} \log \bigl(n e^{\beta \max_i x_i}\bigr) \nonumber \\
&= \beta^{-1}\log n + \max_i x_i. \label{basic}
\end{align}
In other words, $\max x_i \le F_\beta(\bx)\le \beta^{-1}\log n + \max x_i$. 

Thus, taking $\beta \ra \infty$ in (\ref{ineq1}),  we get the
second assertion of the theorem. For the first, note 
that with $\gamma = \max_{1\le i,j\le n} |\gamma^Y_{ij} - \gamma^X_{ij}|$, we have
\[
\biggl|\sum_{1\le i, j\le n} \frac{\partial^2   F_\beta}{\partial
  x_j \partial x_i} (\bx) (\sigma^Y_{ij} -
\sigma^X_{ij})\biggr| \le \frac{\beta \gamma}{2}\sum_{1\le i,j\le
  n}p_i(\bx)p_j(\bx) = \frac{\beta \gamma}{2}.
\]
This shows that
\[
|\ee(F_\beta(\bby)) - \ee(F_\beta(\bbx))| \le \frac{\beta\gamma}{4}.
\]
Combined with (\ref{basic}), this gives
\[
|\ee(\max_i Y_i) - \ee(\max_i X_i)| \le \frac{\beta\gamma}{4} +
\frac{\log n}{\beta}.
\]
Choosing $\beta = 2\sqrt{\frac{\log n}{\gamma}}$ gives the desired
result. 

\begin{small}

\end{small}

\end{document}